\documentclass[reqno,10pt]{amsart}

\usepackage{amsmath,amssymb}
\usepackage[dvips]{graphics}
\usepackage{epsfig}

\newtheorem{theorem}{Theorem}[section]

\newtheorem{condition}{Condition}
\newtheorem{conjecture}[theorem]{Conjecture}
\newtheorem{corollary}[theorem]{Corollary}

\newtheorem{definition}[theorem]{Definition}

\newtheorem{lemma}[theorem]{Lemma}

\numberwithin{equation}{section}

\theoremstyle{definition}
\newtheorem{remark}[theorem]{Remark}

\newcommand\ben{\begin{enumerate}}
\newcommand\een{\end{enumerate}}

\renewcommand{\i}{{\mathrm{i}}} 

\newcommand{\twocase}[5]{#1 \begin{cases} #2 & \text{#3}\\ #4
&\text{#5} \end{cases}   }

\newcommand\be{\begin{equation}}
\newcommand\ee{\end{equation}}
\newcommand\bea{\begin{eqnarray}}
\newcommand\eea{\end{eqnarray}}

\newcommand{\R}{\mathbb{R}}
\newcommand{\C}{\mathbb{C}}
\newcommand{\Z}{\mathbb{Z}}
\newcommand{\Q}{\mathbb{Q}}

\newcommand{\ga}{\alpha}     
\newcommand{\gep}{\epsilon}  

\newcommand{\foh}{\frac{1}{2}}  

\renewcommand{\i}{{\mathrm{i}}} 

\newcommand{\osm}{\overline{S}_m}

\newcommand{\oosm}{\frac1{\sqrt{m}}}
\newcommand{\feta}[1]{\eta\left(#1\right)}

\numberwithin{equation}{section}

\begin{document}

\title[Benford's Law, Values of $L$-Functions and the $3x+1$ Problem]{Benford's Law,
Values of $L$-Functions and the $3x+1$ Problem}
\author{Alex V. Kontorovich}
\email{alexk@math.columbia.edu}
\address{Department of Mathematics, Columbia, New York, NY, $10027$}
\author{Steven J. Miller}
\email{sjmiller@math.ohio-state.edu} \address{Department of
Mathematics, The Ohio State University, Columbus, OH
43210}\curraddr{Department of Mathematics, Brown University,
Providence, RI $02912$} \subjclass[2000]{11K06, 60A10, 11B83,
11M06 (primary), 60F05, 11J86, 60J65, 46F12 (secondary).}
\keywords{Benford's Law, Poisson Summation, $L$-Function, $3x+1$,
Equidistribution, Irrationality Type}
\date{\today}

\begin{abstract}
We show the leading digits of a variety of systems satisfying
certain conditions follow Benford's Law. For each system proving
this involves two main ingredients. One is a structure theorem of
the limiting distribution, specific to the system. The other is a
general technique of applying Poisson Summation to the limiting
distribution. We show the distribution of values of $L$-functions
near the central line and (in some sense) the iterates of the
$3x+1$ Problem are Benford.
\end{abstract}

\maketitle

\section{Introduction}

While looking through tables of logarithms in the late 1800s, Newcomb \cite%
{New} noticed a surprising fact: certain pages were significantly
more worn than others. People were referencing numbers whose
logarithm started with 1 more frequently than other digits. In
1938 Benford \cite{Ben} observed the same digit bias in a wide
variety of phenomenon.

Instead of observing one-ninth (about 11\%) of entries having a
leading digit of 1, as one would expect if the digits $1, 2,
\dots, 9$ were equally likely, over 30\% of the entries had
leading digit 1, and about 70\% had leading digit less than 5.
Since $\log _{10}2\approx 0.301 $ and $\log _{10}5\approx 0.699$,
one may speculate that the probability of observing a digit less
than $k$ is $\log _{10}k$, meaning that the probability of seeing
a particular digit $j$ is $\log _{10}\left( j+1\right) -\log
_{10}j=\log _{10}\left( 1+\frac{1}{j}\right)$. This logarithmic
phenomenon became known as Benford's Law after his paper
containing extensive empirical evidence of this distribution in
diverse data sets gained popularity. See \cite{Hi1} for a
description and history, \cite{Hi2,BBH} for some recent results,
and page 255 of \cite{Knu} for connections between Benford's law
and rounding errors in computer calculations.

In \cite{BBH}  it was proved that many dynamical systems are
Benford, including most power, exponential  and rational
functions, linearly-dominated systems, and non-autonomous
dynamical systems. This adds to the ever-growing family of systems
known or believed to satisfy Benford's Law, such as physical
constants, stock market indices, tax returns, sums and products of
random variables, the factorial function and Fibonacci numbers,
just to name a few.

We introduce two new additions to the family, the Riemann zeta
function (and other $L$-functions) and the $3x+1$ Problem (and
other $(d,g,h)$-Maps), though we prove the theorems in sufficient
generality to include other systems. Roughly, the distribution of
digits of values of $L$-functions near the critical line and the
ratio of observed versus predicted values of iterates of the
$3x+1$ Map tend to Benford's Law. For exact statements of the
results, see Theorem \ref{thm:benfordLfns} and Corollary
\ref{cor:benfordLfns} for $L$-functions and Theorem
\ref{thm:3x+1isBenford} for the $3x+1$ Problem. While the best
error terms just miss proving Benford behavior for $L$-functions
on the critical line, we show that the values of the
characteristic polynomials of unitary matrices are Benford in
Appendix \ref{sec:appvalcharpoly}; as these characteristic
polynomials are believed to model the values of $L$-function, this
and our theoretical results naturally lead to the conjecture that
values of $L$-functions on the critical line are Benford.

A standard method of proving Benford behavior is to show the
logarithms of the values become equidistributed modulo 1; Benford
behavior then follows by exponentiation. There are two needed
inputs. For both systems the main term of the distribution of the
logarithms is a Gaussian, which can be shown to be equidistributed
modulo $1$ by Poisson summation. The second ingredient is to
control the errors in the convergence of the distribution of the
logarithms to Gaussians. For $L$-functions this is accomplished by
Hejhal's refinement of the error terms (his result follows from an
analysis of high moments of integrals of $\log |L(s,f)|$), and for
the $3x+1$ Problem it involves an analysis of the discrepancy of
the sequence $k \log_B 2 \bmod 1$ (which follows from $\log_B 2$
is of finite type; see below).

The reader should be aware that the standard notations from number
theory and probability theory sometimes conflict (for example,
$\sigma$ is used to denote the real part of a point in the complex
plane as well as the standard deviation of a distribution); we try
and follow common custom as much as possible. We denote the
Fourier transform (or characteristic function) of $f$ by
$\widehat{f}\left( y\right) =\int_{-\infty }^{\infty }f\left(
x\right) e^{-2\pi \i xy}dx$. Recall $g(T) = o(1)$ means $g(T)\to
0$ as $T\to \infty$, and $g(T) \ll h(T)$ or $g(T) = O(h(T))$ means
there is some constant $C$ such that for all $T$ sufficiently
large, $|g(T)| \le C h(T)$. Our proof of the Benford behavior of
the $3x+1$ problem uses the (irrationality) type of $\log_B 2$ to
control the errors; a number $\ga$ is of type $\kappa$ if $\kappa$
is the supremum of all $\gamma$ with
\be\label{eq:irrtypedef}\underline{\lim}_{q\to\infty} q^{\gamma +
1} \min_p \left| \ga - \frac{p}{q}\right|\ =\ 0.\ee By Roth's
theorem, every algebraic irrational is of type $1$. See for
example \cite{HS,Ro} for more details.

\section{Benford's Law}

To study leading digits, we use the mantissa function, a
generalization of scientific notation. Fix a base $B>1$ and for a
real number $x>0$ define the mantissa function, $M_{B}\left(
x\right)$, from the unique representation of $x$ by
\begin{equation}
x\ =\ M_B\left( x\right) \cdot B^{k},\text{ with }k\in
\mathbb{Z}\text{ and }M_B\left( x\right) \in \left[ 1,B\right) .
\end{equation}
We extend the domain of mantissa to all of $\mathbb{C}$ via
\begin{equation}
\twocase{M_B\left( x\right) \ =\ }{0}{if $x=0$}{ M_B(|x|)}{if $x
\neq 0$.}
\end{equation}

We study the mantissa of many different types of processes
(discrete, continuous  and mixed), and it is convenient to be able
to use the same language for all.  Take an ordered total space
$\Omega$, for example $\mathbb{N}$ or $ \mathbb{R}^{+}$, and a
(weak notion of)\ measure $\mu $ on $\Omega $  such as the
counting measure or Lebesgue measure. For a subset $A\subset
\Omega $ and an element $T\in \Omega $, denote by $ A_{T}=\left\{
\omega \in A:\omega \leq T\right\} $ the truncated set. We define
the probability of $A$ via density in $\Omega $:

\begin{definition}\label{defi:prob}\label{defi:contben}
$\mathbb{P}\left( A\right) =\lim\limits_{T\rightarrow \infty
}\frac{\mu \left( A_{T}\right) }{\mu \left( \Omega _{T}\right) }$,
provided the limit exists.
\end{definition}

For $A\subset \mathbb{N}$ and $\mu $ the counting measure,
$\mathbb{P}\left( A\right) =\lim\limits_{T\rightarrow \infty
}\frac{\#\left\{ n\in A:\ n\leq T\right\} }{T}$, while if
$A\subset \mathbb{R}^{+}$ and $\mu $ is Lebesgue measure then
$\mathbb{P}\left( A\right) =\lim\limits_{T\rightarrow \infty }
\frac{\mu \left( 0\leq t\leq T:\ t\in A\right) }{T}$. In Appendix
\ref{sec:appvalcharpoly}  we extend our notion of probability to a
slightly more general setting, but this will do for now.

For a sequence of real numbers indexed by $\Omega$,
$\overrightarrow{X}=\left\{ x_{\omega }\right\}_{\omega \in \Omega
}$, and a fixed $s\in \left[ 1,B\right)$, consider the pre-image
of mantissa, $\{ \omega \in \Omega :1\leq M_B( x_{\omega }) \leq
s\} $; we abbreviate this by $ \{ 1\leq M_B( \overrightarrow{X})
\leq s\}$.

\begin{definition}\label{defi:discben}
A sequence $\overrightarrow{X}$ is said to be \textbf{Benford}
(base $B$) if
for all $s\in \left[ 1,B\right) $,%
\begin{equation}
\mathbb{P}\left\{ 1\leq M_B( \overrightarrow{X}) \leq s\right\} \
=\ \log_{B}s.
\end{equation}
\end{definition}

Definition \ref{defi:discben} is applicable to the values of a
function $f$, and we say $f$ is Benford base $B$ if \be
\lim\limits_{T\rightarrow \infty }\frac{ \mu \left( 0\leq t\leq
T:1\leq M_B\left( f\left( t\right) \right) \leq\ s\right)
}{T}=\log_{B}s.\ee

We describe an equivalent condition for Benford behavior which is
based on equidistribution. Recall

\begin{definition} A set $A\subset\R$ is equidistributed modulo $1$ if for any
$[a,b] \subset [0,1]$ we have \be \lim_{T\to \infty }\frac{\mu
\left( \left\{x \in A_T: x \bmod 1 \in [a,b]\right\}\right) }{\mu
\left( A_{T}\right) } \ = \ b-a. \ee \end{definition}

The following two statements are immediate:

\begin{lemma}\label{Lemmakey}
$u\equiv v \bmod 1$ if and only if the mantissa of $B^{u}$ and
$B^{v}$ are the same, base $B$.
\end{lemma}

\begin{lemma}\label{lemIjb}
$y  \bmod 1\in \lbrack 0,\log_{B}s]$ if and only if $B^{y}$ has
mantissa in $ [1,s]$.
\end{lemma}

The following result is a standard way to prove Benford behavior:

\begin{theorem}\label{thm:equidist}
Let $\overrightarrow{Y_B} =\log_{B}| \overrightarrow{X}| $, so
pointwise $y_{\omega,B}=\log_{B}| x_{\omega }| $, and set
$\log_{B}0=0$. Then $\overrightarrow{Y_B}$ is equidistributed
modulo 1 if and only if $\overrightarrow{X}$ is Benford base $B$.
\end{theorem}

\begin{proof}
By Lemma \ref{lemIjb}, the set $\{ \overrightarrow{Y_B}\bmod 1\in %
[ 0,\log_{B}s]\} $ is the same as the set $\{ M_B(
\overrightarrow{X}) \in [ 1,s] \} $. Hence $\overrightarrow{Y_B}$
is equidistributed modulo 1 if and only if
\begin{equation}
\log_{B}s\ =\ \mathbb{P}\left\{ \overrightarrow{Y_B}\bmod 1\in [
0,\log_{B}s] \right\} =\ \mathbb{P}\left\{ M_B(
\overrightarrow{X}) \in \left[ 1,s\right] \right\}
\end{equation}
if and only if $\overrightarrow{X}$ is Benford base $B$.
\end{proof}

Theorem \ref{thm:equidist} reduces investigations of Benford's Law
to equidistribution modulo $1$, which we analyze below.

\begin{remark} The limit in Definition \ref{defi:contben}, often
called the natural density, will exist for the sets in which we
are interested, but need not exist in general. For example, if $A$
is the set of positive integers with first digit $1$, then
$\frac{\#\left\{ n\in A:\ n\leq T\right\} }{T}$ oscillates between
its $\liminf$ of $\frac19$ and its $\limsup$ of $\frac59$. One can
study such sets by using instead the analytic density \be
\mathbb{P}_{\text{an}}\left( A\right) \ = \ \lim_{s\to 1^+} \frac{
\sum_{n\in A} n^{-s} }{\zeta(s)}, \ee where $\zeta(s)$ is the
Riemann Zeta Function (see \S\ref{sec:valuesLFns}). A
straightforward argument using analytic density gives Benford-type
probabilities. In particular, Bombieri (see \cite{Se}, page 76)
has noted that the analytic density of primes with first digit $1$
is $\log_{10}2$, and this can easily be generalized to Benford
behavior for any first digit.
\end{remark}

\section{Poisson Summation and Equidistribution modulo
$1$}\label{sec:poisssum}

We investigate systems $\overrightarrow{X_{T}}$ converging to a
system $\overrightarrow{X}$  with associated logarithmic processes
$\overrightarrow{Y_{T,B}}$. For example, take some function $g:\R
\to \C$ and let $\overrightarrow{X} = \{g(t)\}_{t\in\R}$. Then
$\overrightarrow{X_T} = \left\{g(t) \right\}_{0\leq t\leq T}$ are
truncations of $\overrightarrow{X}$, with log-process
$\overrightarrow{Y_{T,B}} = \left\{\log_B |g(t)| \right\}_{0\leq
t\leq T}$. When there is no ambiguity we drop the dependence on
$B$ and write just $\overrightarrow{Y_{T}}$ for
$\overrightarrow{Y_{T,B}}$.

Let $f(x)$ be a fixed probability density with cumulative
distribution function $F\left( x\right) =\int_{-\infty
}^{x}f\left( t\right) dt$. In our applications the probability
densities of $\overrightarrow{Y_{T,B}}$ are approximately a spread
version of $f$ such as $f_T(x) = \frac{1}{T}f\left( \frac{x}{T}
\right)$. There is, however, an error term, and the log-process
$\overrightarrow{Y_{T,B}}$ has a cumulative distribution function
given by
\begin{eqnarray}\label{eq:proddistEYT}
F_{T}\left( x\right) & \ = \ & \mathbb{P}\left\{
\overrightarrow{Y_{T,B}} \leq x \right\} \nonumber\\
& = & \int_{-\infty }^{x}\frac{1}{T}f\left( \frac{t}{T}
\right) dt+E_{T}\left( x\right) \nonumber\\
&=&F\left( \frac{x}{T}\right) +E_{T}\left( x\right) ,
\end{eqnarray}
where $E_{T}$ is an error term. Our goal is to show that, under
certain conditions, the error term is negligible and $f_{T}\left(
x\right)$ spreads to make $\overrightarrow{Y_{T,B}}$
equidistributed modulo $1$ as $T\rightarrow \infty$. This will
imply that $\overrightarrow{X}$ is Benford base $B$.

In our investigations we need the density $f$, cumulative
distribution function $F_T$ and errors $E_{T}$ to satisfy certain
conditions in order to control the error terms.

\begin{definition}[Benford-good]
Systems $\overrightarrow{Y_{T,B}}$ with cumulative distribution
functions $F_T$ are \textbf{Benford-good} if the $F_T$ satisfy
\eqref{eq:proddistEYT}, the probability density $f$ satisfies
sufficient conditions for Poisson Summation ($\sum_n f(n) = \sum_n
\widehat{f}(n)$), and there is a monotone increasing function
$h(T)$ with $\lim_{T\to\infty} h(T) = \infty$ such that $f$ and
$E_T$ satisfy

\begin{condition}\label{cond:1}  Small tails:
\be F_T(\infty) - F_T(Th(T)) \ = \ o(1),\ \  \ \ \ F_T(-Th(T)) -
F_T(-\infty) \ = \ o(1). \ee
\end{condition}

\begin{condition}\label{cond:3} Rapid decay of the characteristic function:
\be S\left( T\right)\ = \ \sum_{k\neq 0}\left|
\frac{\widehat{f}(Tk)}{k}\right| \ = \ o(1). \ee
\end{condition}

\begin{condition}\label{cond:4} Small truncated translated error:
\be \mathcal{E}_{T}(a,b)\ = \ \sum\limits_{|k| \le
Th(T)}\left[E_{T}( b+k) - E_T(a+k)\right] \ = \ o(1), \ee for all
$0\le a < b \le 1$.
\end{condition}
\end{definition}

In all our applications $f$ will be a Gaussian, in which case the
Poisson Summation Formula holds. See for example \cite{Da} (pages
14 and 63).

Condition \ref{cond:1} asserts that essentially all of the mass
lies in $[-Th(T),Th(T)]$. In applications $T$ will be the standard
deviation, and this will follow from Central Limit type
convergence.

Condition \ref{cond:3} is quite weak, and is satisfied in most
cases of interest. For example, if $f$ is differentiable and
$f^{\prime }$ is integrable\ (as is the case if $f$ is the
Gaussian density), then $|\widehat{f}(y)|\leq \frac{1}{|y|} \int
|f^{\prime }(x)|dx=O\left( \frac{1}{|y|}\right) $, which suffices
to show $S\left( T\right) = o(1)$.

Condition \ref{cond:4} is the most difficult to prove for a
system, and to our knowledge has not previously been analyzed in
full detail. It is well known (see \cite{Fe}) that there are some
processes (for example, Bernoulli trials) with standard deviation
of size $T$ where the best attainable estimate is $E_T(x) =
O\left(\frac1{T}\right)$. Errors this large lead to
$\mathcal{E}_T(a,b) = O(1)$.

We now see why these conditions suffice. For $[a,b]\subset \lbrack
0,1)$, let $P_T(a,b)$ denote the probability that
$\overrightarrow{Y_{T,B}}\bmod 1\in \left[ a,b\right]$. To prove
$\overrightarrow{Y_{T,B}}$ becomes equidistributed modulo $1$, we
must show that $P_T[a,b] \rightarrow b-a$. We would like to argue
as follows:
\begin{eqnarray}\label{eq:prepoisswitherr}
P_T[a,b]  &\ = \ &\mathbb{P}\left\{ \overrightarrow{Y_{T,B}}\bmod
1\in \left[ a,b\right] \right\}\nonumber\\ & \ = \ & \sum_{k\in
\mathbb{Z}}\mathbb{P}\left\{
\overrightarrow{Y_{T,B}}\in \left[ a+k,b+k\right] \right\} \  \nonumber\\
&\ = \ &\sum_{k\in \mathbb{Z}}\left( F_{T}\left( b+k\right)
-F_{T}\left(
a+k\right) \right) \nonumber\\
&\ = \ &\sum_{k\in \mathbb{Z}}\left[\int_{a}^{b}\frac{1}{T}f\left(
\frac{x+k}{T} \right) dx+E_{T}( b+k)
-E_{T}( a+k)\right] \nonumber\\
&\ = \ &\sum_{k\in \mathbb{Z}}\left[\int_{a}^{b}\frac{1}{T}f\left(
\frac{x+k}{T} \right) dx\right]+\sum_{k\in \mathbb{Z}}
\left[E_{T}(b+k) -E_{T}(a+k)\right].
\end{eqnarray}

While the main term can be handled by a straightforward
application of Poisson Summation, the best pointwise bounds for
the error term are not summable over all $k\in\Z$. This is why
Condition \ref{cond:1} is necessary, so that we may restrict the
summation.

\begin{theorem}\label{thm:poissum} Assume log-processes $\overrightarrow{Y_{T,B}}$
are \textbf{Benford-good}. Then
$\overrightarrow{Y_{T,B}}\rightarrow \overrightarrow{Y_B}$, where
$\overrightarrow{Y_B}$ is equidistributed modulo 1.
\end{theorem}

\begin{proof}
As the Fourier transform converts translation to multiplication,
if $ g_{x}(u)=f\left( \frac{u+x}{T}\right)$ then a straightforward
calculation shows that $\widehat{g_{x}}(w)=e^{2\pi \i
xw}T\widehat{f}(Tw)$ for any fixed $x$. Our assumptions on $f$
allow us to apply Poisson Summation to $g$, and we find
\begin{equation}\label{eq:Poissum}
\sum_{k\in \mathbb{Z}} f \left( \frac{x+k}{T}\right)\ =\
\sum_{k\in \mathbb{Z}}g_{x}(k)\ =\ \sum_{k\in
\mathbb{Z}}\widehat{g_{x}}(k)\ =\ T\sum_{k\in \mathbb{Z}}e^{2\pi
\i xk}\widehat{f}(Tk).
\end{equation}
Let $[a,b] \subset [0,1]$. By Condition \ref{cond:1} and
\eqref{eq:proddistEYT}, \bea P_T(a,b)& \ = \ & \sum_{|k| \le
Th(T)} \left( F_T(b+k) - F_T(a+k) \right)\nonumber\\ & & \ \ \ \ \
\ \ \ \ \ +\ O\left( F_T(\infty) - F_T(Th(T)) \right) + O\left(
F_T(-Th(T)) - F_T(-\infty) \right) \nonumber\\ & = & \sum_{|k| \le
Th(T)}
 \left[ \frac1{T}\int_a^b f\left(\frac{x+k}{T}\right)dx +
 E_T(b+k) - E_T(a+k)  \right] + o(1) \nonumber\\ & = & \sum_{|k| \le Th(T)}
\frac1{T}\int_a^b f\left(\frac{x+k}{T}\right)dx +
 \mathcal{E}_T(a,b)   + o(1). \eea

By Condition \ref{cond:4}, $\mathcal{E}_{T}(a,b) = o(1)$; as $f$
is integrable we may return the sum to all $k\in \Z$ at a cost of
$o(1)$. The interchange of summation and integration below is
justified from the decay properties of $f$. To see this, simply
insert absolute values in the arguments. Therefore using
\eqref{eq:Poissum},
\begin{eqnarray}\label{eq:proofaftererrafterpoiss}
P_T[a,b] &\ = \ &\frac{1}{T}\sum_{k\in \mathbb{Z}
}\int_{a}^{b}f\left( \frac{x+k}{T}\right) dx + o(1)\nonumber\\ & =
& \frac{1}{T}\int_{a}^{b}\left(
\sum_{k\in \mathbb{Z}}g_{x}\left( k\right) \right) dx + o(1) \nonumber\\
&=&\frac{1}{T}\int_{a}^{b}\left( \sum_{k\in
\mathbb{Z}}\widehat{g_{x}} \left( k\right) \right) dx  + o(1)\nonumber\\
& =  & \sum_{k\in \mathbb{Z}}\widehat{f}
(Tk)\int_{a}^{b}e^{2\pi \i xk}dx  + o(1)\nonumber\\
&=&\widehat{f}\left( 0\right) \left( b-a\right) +\sum_{k\neq
0}\widehat{f} \left( Tk\right) \frac{e^{2\pi \i bk}-e^{2\pi \i
ak}}{2\pi \i k}  + o(1).
\end{eqnarray}
As $f$ is a probability density, $\widehat{f}(0)=1$, and by
Condition \ref{cond:3} the sum in
\eqref{eq:proofaftererrafterpoiss} is $o(1)$. Therefore \be
P_T(a,b) \ = \ b-a + o(1), \ee which completes the proof.
\end{proof}

As an immediate consequence, we have:

\begin{theorem}\label{thm:main}
Let $\overrightarrow{X_{T}}$ (the truncation of
$\overrightarrow{X}$) have corresponding log-process
$\overrightarrow{Y_{T,B}}$. Assume the $\overrightarrow{Y_{T,B}}$
are \textbf{Benford-good}. Then $\overrightarrow{X}$ is Benford
base $B$.
\end{theorem}

\begin{proof} This follows immediately from Theorems
\ref{thm:poissum} and \ref{thm:equidist}. \end{proof}

An immediate application of Theorem \ref{thm:main} is to processes
where the distribution of the logarithms is exactly a spreading
Gaussian (i.e., there are no errors to sum). We describe such a
situation below.

Recall a \emph{Brownian motion} (or Wiener process) is a
continuous process with independent, normally distributed
increments. So if $W$ is a Brownian motion, then $W_{t} - W_{s}$
is a random variable having the Gaussian distribution with mean
zero and variance $t-s$, and is independent of the random variable
$W_{s} - W_{u}$ provided $u<s<t$.

A standard realization of Brownian motion is as the scaled limit
of a random walk. Let $x_{1}, x_{2}, x_{3}, \dots$ be independent
Bernoulli trials (taking the values $+1$ and $-1$ with equal
probability) and let $S_{n}=\sum_{i=1}^{n} x_{i}$ denote the
partial sum. Then the normalized process \be W_{t}^{(n)} \ =\
\frac{1}{\sqrt{n}}\ S_{nt} \ee (extended to a continuous process
by linear interpolation) converges as $n \to \infty$ to the Wiener
process. See \cite{Bi} or Chapter 2.4 of \cite{KaSh} for further
details.

A \emph{geometric Brownian motion} is simply a process $Y$ such
that the process $\log Y$ is a Brownian motion. It was known to
Benford that stock market indices empirically demonstrated this
digit bias, and for almost as long these indices have been
modelled by geometric Brownian motion. Thus Theorem \ref{thm:main}
implies the well-known observation that

\begin{corollary}
A geometric Brownian motion is Benford.
\end{corollary}

\section{Values of $L$-Functions}\label{sec:valuesLFns}

Consider the Riemann Zeta function
\begin{equation}
\zeta \left( s\right)\ =\ \sum_{n=1}^{\infty }\frac{1}{n^{s}} \ =
\ \prod_{p \ \text{prime}} \left(1 - \frac1{p^s}\right)^{-1}.
\end{equation} Initially defined for ${\rm Re}(s) > 1$, $\zeta(s)$
has a meromorphic continuation to all of $\C$. More generally, one
can study an $L$-function  \be L\left( s,f\right) \ = \
\sum_{n=1}^{\infty }\frac{a_{f}(n)}{n^{s}} \ = \  \prod_{p\ {\rm
prime}}\ \prod_{j=1}^d\left( 1-\frac{\alpha
_{f,d}(p)}{p^{s}}\right)^{-1},
\end{equation} where the coefficients $a_f(n)$ have arithmetic
significance. Common examples include Dirichlet $L$-functions
(where $a_f(n) = \chi(n)$ for a Dirichlet character $\chi$) and
elliptic curve $L$-functions (where $a_f(p)$ is related to the
number of points on the elliptic curve modulo $p$).

All the $L$-functions we study satisfy (after suitable
renormalization) a functional equation relating their value at $s$
to their value at $1-s$. The region $0 \le {\rm Re}(s) \le 1$ is
called the critical strip, and ${\rm Re}(s) = \foh$ the critical
line. The behavior of $L$-functions in the critical strip,
especially on the critical line, is of great interest in number
theory. The Generalized (or, as some prefer, Grand) Riemann
Hypothesis, GRH, asserts that the zeros of any ``nice''
$L$-function are on the critical line. The location of the zeros
of $\zeta(s)$ is intimately connected with the error estimates in
the Prime Number Theorem. The Riemann Zeta function can be
expressed as the moment of the maximum of a Brownian Excursion,
and the distribution of the zeros (respectively, values) of
$L$-functions is believed to be connected to that of eigenvalues
(respectively, values of characteristic polynomials) of random
matrix ensembles. See \cite{BPY,Con,KaSa,KeSn} for excellent
surveys.

We investigate the leading digits of $L$-functions near the
critical line, and show that the distribution of the digits of
their absolute values is Benford (see Theorem
\ref{thm:benfordLfns} for the precise statement). The starting
point of our investigations of values of the Riemann zeta function
along the critical line $s=\tfrac{1}{2}+it$ is the log-normal law
(see \cite{Lau, Sel1}):
\begin{equation}
\lim_{T\rightarrow \infty }\frac{\mu \left( \left\{ 0\leq t\leq
T:\log
|\zeta \left( \tfrac{1}{2}+\i t\right)| \ \leq \ y\sqrt{%
\tfrac{1}{2}\log \log T}\right\} \right) }{T}\ =\
\frac{1}{\sqrt{2\pi}} \int_{-\infty }^{y}e^{-u^{2}/2}du.
\end{equation}
Thus the density of values of $\log \left\vert \zeta \left(
\tfrac{1}{2} +\i t\right) \right\vert $ for $t\in \lbrack 0,T]$
are well approximated by a Gaussian with mean zero and standard
deviation \be \psi_T\ = \ \sqrt{\tfrac{1}{2}\log \log T}+O(\log
\log \log T).\ee Such results are often used to investigate small
values of $|\zeta \left( \frac{1}{2}+\i t\right)|$ and gaps
between zeros. As such, the known error terms are too crude for
our purposes. In particular, one has (trivially modifying (4.21)
of \cite{Hej} or (8) of \cite{Iv}) that
\begin{equation}\label{eq:normalHejIv}
\frac{\mu \left( \left\{ t\in \lbrack T,2T]:a\leq \log \left\vert
\zeta \left( \tfrac{1}{2}+\i t\right) \right\vert \leq b\right\}
\right) }{T} = \frac{1}{\sqrt{2\pi
\psi_{T}^{2}}}\int_{a}^{b}e^{-u^{2}/2\psi_{T}^{2}}du+O\left(
\frac{\log ^{2}\psi_{T}}{\psi_{T}}\right).
\end{equation}
The main term is Gaussian with increasing variance, precisely what
we require for equidistribution modulo 1. The error term, however,
is too large for pointwise evaluation (as we have of the order
$\psi_{T}\log \psi_{T}$ intervals $ [a+n,b+n]$).

Better pointwise error estimates are obtained for many
$L$-functions in \cite{Hej}. These estimates are good enough for
us to see Benford behavior as $T\rightarrow \infty $ near the line
$\text{\rm Re}(s)= \frac{1}{2}$. Explicitly, consider an
$L$-function (or a linear combination of $L$-functions, though for
simplicity of exposition we confine ourselves to the case of one
$L$-function) satisfying

\begin{definition}[Good $L$-Function]\label{defi:goodLfn}
We say an $L$-function is good if it satisfies the following
properties:
\begin{enumerate}
\item Euler product: \be\label{eq:hejcond1} L(s,f)\ = \
\sum_{n=1}^{\infty } \frac{a_{f}(n)}{n^{s}}\ = \ \prod_{p\ {\rm
prime}} \ \prod_{j=1}^d \left(
1-\alpha_{f,j}(p)p^{-s}\right)^{-1}.\ee

\item $L(s,f)$ has a meromorphic continuation to $\mathbb{C}$, is
of finite order, and has at most finitely many poles (all on the
line ${\rm Re}(s)=1$).

\item Functional equation: \be e^{\i\omega }G(s)L(s,f)\ = \
e^{-\i\omega }\overline{G(1-\overline{s})L(1-%
\overline{s})},\ee where $\omega \in \mathbb{R}$ and \be G(s)\ = \
Q^{s} \prod_{i=1}^{h}\Gamma (\lambda _{i}s+\mu _{i})\ee with
$Q,\lambda _{i}>0$ and ${\rm Re}(\mu _{i})\geq 0$.

\item For some $\aleph >0$, $c\in \mathbb{C}$, $x\geq 2$ we have
\be \sum_{p\leq x} \frac{|a_{f}(p)|^{2}}{p}\ = \ \aleph \log \log
x+c+O\left( \frac1{\log x}\right).\ee

\item The $\alpha _{f,j}(p)$ are (Ramanujan-Petersson) tempered: $
|\alpha _{f,j}(p)|\leq 1$.

\item If $N(\sigma ,T)$ is the number of zeros $\rho $ of $L(s)$
with ${\rm Re}(\rho )\geq \sigma $ and ${\rm Im} (\rho )\in
\lbrack 0,T]$, then for some $\beta > 0$ we have
\be\label{eq:hejcond6} N(\sigma ,T)\ = \ O\left( T^{1-\beta \left(
\sigma -\tfrac{1}{2}\right) }\log T\right).\ee
\end{enumerate} \end{definition}

\begin{remark}\label{rek:goodLfns}
There are many families of $L$-functions which satisfy the above
six conditions. The last two are the most difficult conditions to
verify, as in all cases where these are known the first four
conditions can be shown to be satisfied. The last two conditions
are established for many $L$-functions (for example, see
\cite{Sel1} for $\zeta(s)$ and \cite{Luo} for holomorphic Hecke
cuspidal forms of full level and even weight $k>0$; see Chapter 10
\cite{IK} for more on the subject), and is an immediate
consequence of GRH.
\end{remark}

We quote a version of the log-normal law with better error terms
(see (4.20) from \cite{Hej} with a trivial change of variables in
the Gaussian integral); for the convenience of the reader we list
where the various parameters in Hejhal's result are defined. The
error terms will be pointwise summable, and allow us to prove
Benford behavior.

\begin{theorem}[Hejhal]\label{thm:HejLogNorm}
Let $L(s,f)$ be a good $L$-function as in Definition
\ref{defi:goodLfn}, and

\begin{itemize} \item fix $\delta \in (0,1)$ (\cite{Hej},
Lemmas 2 and 3, page 556), $g \in (0,1] $ (\cite{Hej}, Lemma 3,
page 556) and $\kappa \in (1,3]$ (\cite{Hej}, page 560 and (4.18)
on page 562);

\item choose $\sigma \geq \tfrac{1}{2}+\frac{g}{\log y}$
(\cite{Hej}, page 563) and
$\tfrac{1}{2}\leq \sigma \leq \tfrac{1}{2}+\frac{1}{\log ^{\delta }T}$ (\cite%
{Hej}, page 562);

\item the variance $\psi (\sigma ,T)$ (see \cite{Hej}, Lemma 1,
page 566) satisfies \be\psi (\sigma ,T)\ = \ \aleph \log \left[ \min \left( \log T,\frac{1}{%
\sigma -\tfrac{1}{2}}\right) \right] +O(1); \ee

\item choose $N=\lfloor \psi (\sigma ,T)^{\kappa }\rfloor $ and $y=T^{1/2N}$ (\cite%
{Hej}, (4.18), page 565).
\end{itemize} Then we have
 \bea\label{eq:lognormalhej} & & \frac{\mu\left( \left\{t \in
[T,2T]: a \le \log \left|L\left(\sigma + \i t,f\right)\right| \le
b \right\} \right)}{T} \ = \ \frac{1}{\sqrt{\psi(\sigma, T)}}
\int_a^b e^{-\pi u^2/\psi(\sigma, T)}du\nonumber\\ & & \ \ \ + \
O\left( \frac{1}{\psi(\sigma, T)} \min\left(1,
\frac{|b-a|}{\sqrt{\psi(\sigma,T)}}\right) +
\psi(\sigma,T)^{-\kappa/2} + y^{(1/3)(1-2\sigma)}\right), \eea the
implied constant depends only on $\beta$ (Condition (6) of
Definition \ref{defi:goodLfn}), $f$, $\delta$, $g$ and $\kappa$.
 \end{theorem}

For our purposes, a satisfactory choice is to take $\sigma
=\tfrac{1}{2}+\frac{1}{\log ^{\delta }T}$ and $\kappa >2$. Then
$\psi (\sigma ,T)=\aleph \log \log
T+O(1) $ and%
\begin{eqnarray}
y^{(1/3)(1-2\sigma )}\ &\ = \ &\ T^{\frac{1}{\log ^{\delta
}T}\frac{-1}{3(\aleph \log \log T+O(1))^{\kappa }}}\ = \ \exp
\left(
-\frac{\log ^{1-\delta }T}{3(\aleph \log \log T+O(1))^{\kappa }}\right) \nonumber\\
&\ll &\frac{\left( \log \log T\right) ^{\kappa }}{\log ^{1-\delta
}T}.
\end{eqnarray}

We now show, in a certain sense, the values of $|L(s,f)|$ are
Benford. While any modest cancellation would yield the following
result on the critical line, due to our error terms for each
interval $[T,2T]$ we must stay slightly to the right of ${\rm
Re}(s)=\frac{1}{2}$.

\begin{theorem}\label{thm:benfordLfns} Let $L(s,f)$ be a good
$L$-function as in Definition \ref{defi:goodLfn}; for example we
may take $\zeta(s)$. If the GRH and Ramanujan conjectures hold we
may take any cuspidal automorphic $L$-function; see also Remark
\ref{rek:goodLfns}. Fix a $\delta \in (0,1)$. For each $T$, let
$\sigma_T=\frac{1}{2}+\frac{1}{\log ^{\delta }T}$. Then
\begin{equation}
\lim\limits_{T\rightarrow \infty }\frac{\mu \left\{ t\in \lbrack
T,2T]:\ 1\leq M_B\left( |L(\sigma_T +\i t,f)|\right) \leq \tau
\right\} }{T}\ =\ \log_{B}\tau .
\end{equation}
Thus the values of the $L$-function satisfy Benford's Law in the
limit (with the limit taken as described above) for any base $B$.
\end{theorem}

\begin{proof} We first prove the claim for base $e$, and then comment on the
changes needed for a general base $B$. Unfortunately the notation
from number theory slightly conflicts with the standard notation
from probability theory of \S\ref{sec:poisssum}. By Theorem
\ref{thm:equidist}, it suffices to show that \be
\lim\limits_{T\rightarrow \infty }\frac{\mu \left\{ t\in \lbrack
T,2T]:\ a \leq \log |L(\sigma_T +\i t,f)| \leq b \right\} }{T}\ =\
b - a. \ee

Let $\psi_T =\psi(\sigma_T,T)$ be the variance of the Gaussian in
\eqref{eq:lognormalhej}, which tends to infinity with $T$. The
standard deviation is thus $\sqrt{\psi_T}$, and corresponds to
what we called $T$ in \S\ref{sec:poisssum}. Let $\eta(x)$ be the
standard normal (mean zero, variance one; $\eta$ plays the role of
$f$ from \S\ref{sec:poisssum} -- as it is standard to denote
$L$-functions by $L(s,f)$, we use $\eta$ here and in
\S\ref{sec:3x+1}), and set $\eta_{\sqrt{\psi_T}}(x) =
\frac1{{\sqrt{\psi_T}}}\
\eta\left(\frac{x}{{\sqrt{\psi_T}}}\right)$. Note
$\eta_{\sqrt{\psi_T}}(x)$ is the density of a normal with mean
zero and variance $\psi_T$. By \eqref{eq:lognormalhej} we have \be
F_T(x) \ = \ \int_{-\infty}^x \eta_{\sqrt{\psi_T}}(x)dx + E_T(x),
\ee where $E_T(x) = O(\psi_T^{-1})$. We must show the logarithms
of the absolute values of the $L$-function are Benford-good. As
$\eta$ is a Gaussian it satisfies the conditions for the Poisson
Summation Formula, and the log-process $\overrightarrow{Y_T} =
\log\left|L(\sigma_T+\i t, f)\right|$ satisfies
\eqref{eq:proddistEYT}. Thus to apply Theorem \ref{thm:main} it
suffices to show $\eta$, $F_T$ and $\mathcal{E}_T$ satisfy
Conditions \ref{cond:1} through \ref{cond:4} for some monotone
increasing function $h(\psi_T)$ with $\lim_{T\to\infty} h(\psi_T)
= \infty$. We take $h(\psi_T) = \sqrt{\log \psi_T}$.

Condition \ref{cond:1} is immediately verified. To show
$F_{\sqrt{\psi_T}}(\infty) - F_{\sqrt{\psi_T}}({\sqrt{\psi_T}}
h(\psi_T)) = o(1)$ we use \eqref{eq:lognormalhej} to conclude the
contribution from the error is $o(1)$, and then note that the
integral of the Gaussian with standard deviation $\sqrt{\psi_T}$
past $\sqrt{\psi_T\log\psi_T}$ is small (as $\eta$ is the density
of the standard normal, this integral is dominated by \be
\frac1{\sqrt{2\pi}} \int_{|x| \ge \sqrt{\log \psi_T}}\eta(x)dx,
\ee which is $o(1)$). Identical arguments show
$F_{\sqrt{\psi_T}}(-{\sqrt{\psi_T}} h(\psi_T)) -
F_{\sqrt{\psi_T}}(-\infty) = o(1)$. As we are integrating a
sizable distance past the standard deviation, it is easy to see
that the contribution from the Gaussian is small. We do not need
the full strength of the bounds in \eqref{eq:lognormalhej}; the
bounds from \eqref{eq:normalHejIv} suffice to control the errors.

Condition \ref{cond:3} follows from the trivial fact that $\eta'$
is integrable. We now show Condition \ref{cond:4} holds. Here the
bounds from \eqref{eq:normalHejIv} just fail. Using those bounds
and summing over $|k| \le \sqrt{\psi_T}h(\psi_T)$ would yield an
error of size $O\left(\sqrt{\psi_T} h(\psi_T) \cdot \frac{\log^2
\sqrt{\psi_T}}{\sqrt{\psi_T}}\right) = O\left(\log^{2.5}
\psi_T\right)$. We instead use \eqref{eq:lognormalhej}, and find
for $[a,b] \subset [0,1]$ that \bea \mathcal{E}_T(a,b) & \ = \ &
\sum_{|k|
\le \sqrt{\psi_T} h(\psi_T)}\left[ E_T(b+k) - E_T(a+k)\right] \nonumber\\
& = & \sum_{|k| \le \sqrt{\psi_T\log \psi_T}} O\left(
\frac{1}{\psi_T} \min\left(1,
\frac{|b-a|}{\sqrt{\psi_T}}\right)+\psi_T^{-\kappa/2}  +
y^{(1/3)(1-2\sigma)}\right) \nonumber\\ & = &
O\left(\frac{\sqrt{\log\psi_T}}{\sqrt{\psi_T}}  + \psi_T^{\foh -
\frac{\kappa}{2}} \sqrt{\log\psi_T}  + \sqrt{\psi_T\log
\psi_T}\frac{\left( \log \log T\right) ^{\kappa }}{\log ^{1-\delta
}T} \right) \nonumber\\ & \ =\ &  o(1) \eea because $\kappa > 1$,
$\delta < 1$ and $\psi_T \ll \log\log T$.

As all the conditions of Theorem \ref{thm:poissum} are satisfied,
we can conclude that \be P_{\sqrt{\psi_T}}(a,b) \ = \ b - a +
o(1). \ee  We have shown that tending to infinity in this manner,
the distribution corresponding to $\log | L(\sigma_T +\i t,f)| $
converges to being equidistributed modulo $1$, which by Theorem
\ref{thm:main} implies the values of $|L(\sigma_T +\i t,f)|$ are
Benford base $e$ (as always, along the specified path converging
to the critical line).

For a general base $B$, note $\log_B x = \frac{\log x}{\log B}$.
The effect of changing base is that $\log_B |L(\sigma_T+\i t,f)|$
converges to a Gaussian with mean zero and variance $\frac{1}{\log
B} \cdot \sqrt{\psi(\sigma_T,T)}$ (instead of mean zero and
variance $\sqrt{\psi(\sigma_T,T)}$). The argument now proceeds as
before.
\end{proof}

\begin{corollary}\label{cor:benfordLfns}
Theorem \ref{thm:benfordLfns} is valid if instead of intervals
$[T,2T]$ we consider intervals $[0,T]$.
\end{corollary}

\begin{proof}
Let $\alpha (T)=(\log \log \log T)^{\log 2}$. We consider the
intervals $I_{0} = [0,T/\alpha (T)]$ and
\begin{equation}
I_{i}\ =\ \left[ 2^{i-1}T/\alpha (T),\ 2^{i}T/\alpha (T)\right] ,\
\ i\in \{1,2,\dots ,\log \log \log \log T\}.
\end{equation}

We may ignore $I_{0}$ as it has length $o(T)$. For each interval
$I_{i}$, $i \geq 1$, we use \eqref{eq:lognormalhej} and argue as
before. We may keep the same values of $\beta ,\delta ,g,\kappa
,\sigma_T$ as before. $T$ and $y$ change, which implies $\psi_T
=\psi(\sigma_T,T)$ changes; however, the leading term of $\psi_T$
is still $\aleph \log \log T$, and $ y^{(1/3)(1-2\sigma )}$ again
leads to negligible contributions. As there are only $\log \log
\log \log T$ intervals, we may safely add all the errors.
\end{proof}

\begin{remark} If we stay a fixed distance off the critical line, we
do not expect Benford behavior. This is because for a fixed
$\sigma > \foh$, for $\zeta(s)$ we have a distribution function
$G_\sigma$ such that \be \lim_{T\to\infty} \frac{\mu\{t \in [0,T]:
\log |\zeta(\sigma + \i t)| \in [a,b] \}}{T} \ = \ \int_a^b
G_\sigma(u)du. \ee Unlike the log-normal law
\eqref{eq:normalHejIv}, where the variance increases with $T$,
note here there is no increasing variance for fixed $\sigma$
(though of course the variance depends on $\sigma$); see
\cite{BJ,JW} for proofs. Thus to see Benford behavior it is
essential that as $T$ increases our distance to the critical line
decreases.
\end{remark}

For investigations on the critical line, one can easily show
Benford's Law holds for a truncation of the series expansion of
$\log |L(\frac{1}{2} +\i t,f)| $, where the truncation depends on
the height $T$. See (4.12) of \cite{Hej} for the relevant version
of the log-normal law (which has a significantly better error term
than \eqref{eq:lognormalhej}). Similarly, one can prove statements
along these lines for the real and imaginary parts of
$L$-functions.

Numerical investigations also support the conjectured Benford
behavior. In Figure 1 we plot the percent of first digits of
$\left|\zeta\left(\foh + \i t\right)\right|$ versus the Benford
probabilities for
$t = \frac{k}{4}$, $k \in \{0,1,\dots,65535\}$, and note the
Benford behavior quickly sets in. Of course, we believe that this
is strong evidence for Benford behavior exactly on the critical
line, but as they stand, our error terms are too big and our
cancellation too small to demonstrate this statement.
\begin{figure}\label{fig:plotzetavsbenf}
\begin{center}
\includegraphics[width=10cm]{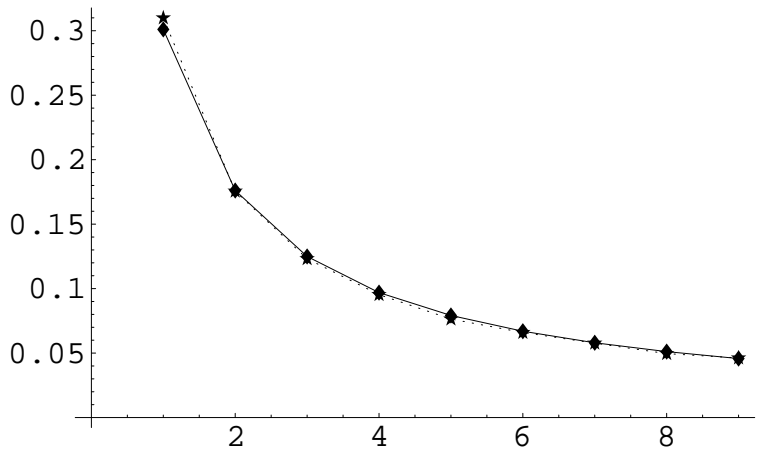}
\caption{Distribution of Digits of $|\zeta(s)|$ versus Benford
Probabilities}
\end{center}\end{figure}

It is believed that values of characteristic polynomials of random
matrix ensembles model values of $L$-functions on the critical
line. In Theorem \ref{thm:valuesCharPoly} of Appendix
\ref{sec:appvalcharpoly} we show that the digit distribution of
the values of these characteristic polynomials converge to the
Benford probabilities (as the size of the matrices tend to
infinity), providing additional support for the conjecture that
$L$-functions are Benford on the critical line.

\section{The $3x+1$ Problem}\label{sec:3x+1}

People working on the
Syracuse-Kakutani-Hasse-Ulam-Hailstorm-Collatz-$\left( 3x+1\right)
$-Problem (there have been a few) often refer to two striking
anecdotes. One is Erd\"{o}s' comment that \textquotedblleft
Mathematics is not yet ready for such problems.\textquotedblright\
The other is Kakutani's communication to Lagarias:
\textquotedblleft For about a month everybody at Yale worked on
it, with no result. A similar phenomenon happened when I\
mentioned it at the University of Chicago. A joke was made that
this problem was part of a conspiracy to slow down mathematical
research in the U.S.\textquotedblright\ Coxeter has offered \$50
for its solution, Erd\"{o}s \$500, and Thwaites, \pounds 1000. The
problem has been connected to holomorphic solutions to functional
equations, a Fatou set having no wandering domain, Diophantine
approximation of $\log _{2}3$, the distribution $\bmod$ $1$ of
$\left\{ \left( \frac{3}{2}\right) ^{k}\right\}_{k=1}^{\infty }$,
ergodic theory on $\mathbb{Z}_{2}$, undecidable algorithms, and
geometric Brownian motion, to name a few (see \cite{Lag1, Lag2}).
We now relate the $\left( 3x+1\right)$-Problem to Benford's Law.

\subsection{The Structure Theorem}

If $x$ is a positive odd integer then $3x+1$ is even, so we can
find an integer $k\geq 1$ such that $2^{k}\parallel \left(
3x+1\right) $, i.e. so that
\begin{equation}
y\ =\ \frac{3x+1}{2^{k}}
\end{equation}
is also odd. In this way, we get the $\left( 3x+1\right) $-Map%
\begin{equation}
M:x\longmapsto y.
\end{equation}
We call the value of $k$ that arises in the definition of $y$  the
\textbf{$k $-value} of $x$. Notice that $y$ is odd and relatively
prime to $3$, so the natural domain for iterating $M$ is the set
$\Pi $ of positive integers prime to $2$ and $3$. Write $ \Pi
=6\mathbb{N}+E$, where $E=\left\{ 1,5\right\} $ is the set of
possible congruence classes modulo 6. The total space is $ \Omega
=\Pi$, not $\mathbb{N}$ or $\mathbb{R}$, and the measure is the
appropriate counting measure.

For every integer $x\in \Pi$ with $0<x<2^{60}$, computers have
verified that enough iterations of the $\left( 3x+1\right) $-Map
eventually send $x$ to the unique fixed point, 1. The natural
conjecture asks if the same statement holds for all $x\in \Pi $:

\begin{conjecture}[$(3x+1)$-Conjecture]
For every $x\in \Pi $, there is an integer $n$ such that
$M^{n}\left( x\right) =1$.
\end{conjecture}

Suppose we apply $M$ a total of $m$ times, calling $x_{0}=x$ and
$x_{i}=M^{i}\left( x\right)$, $i\in\{1,2,\dots,m\}$. For each
$x_{i-1}$ there is a $k$-value, say $k_{i}$, such that
\begin{equation}
x_{i}\ =\ M\left( x_{i-1}\right) \ =\
\frac{3x_{i-1}+1}{2^{k_{i}}}, \ \ \ \ \ i\ \in \ \{1,2,\dots,m\}.
\end{equation}
We store this information in an ordered $m$-tuple $\left(
k_{1},k_{2},\dots,k_{m}\right) $, called the \textbf{$m$-path} of
$x$. Let $\gamma_{m}$ denote the map sending $x$ to its $m$-path,
\begin{equation}
\gamma_{m}:x\mapsto \left( k_{1},k_{2},\dots,k_{m}\right) .
\end{equation}

The natural question is whether given an $m$-tuple of positive
integers $ \left( k_{1},k_{2},\dots,k_{m}\right) $, there is an
integer $x$ whose $m$-path is precisely this $m$-tuple. If so, we
would like to classify the set of all such $x$. In other words, we
want to study the inverse map $\gamma _{m}^{-1}$.

The answer is given by the Structure Theorem, proved in
\cite{KonSi}: for each $m$-tuple $\left(
k_{1},k_{2},\dots,k_{m}\right) $, not only does there exist an $x$
having this $m$-path, but this path is enjoyed by two full
arithmetic progressions, $x\in \left\{
a_{1}n+b_{1},a_{2}n+b_{2}\right\} _{n=0}^{\infty }$, and we can
solve explicitly for $a_{i}$ and $b_{i}$. In fact,
$a_{1}=a_{2}=6\cdot 2^{k_{1}+k_{2}+\cdots+k_{m}}$, and
$b_{i}<a_{i}$ (so the progressions are full; we do not miss any
terms at the beginning). Moreover, the two progressions fall into
the two possible equivalence classes modulo 6; i.e., $\{b_{1}\bmod
6,b_{2}\bmod 6\}=\{1,5\}$. The structure theorem is the key
ingredient in analyzing the limiting distributions. These will
satisfy the conditions of our main theorem (Theorem
\ref{thm:main}), and yield Benford's Law.

Recall (Definition \ref{defi:prob}) that we define the probability
of a subset $A\subset \Pi $ by \bea\label{eq:densDef}
\mathbb{P}\left( A\right) \ =\ \lim_{T\rightarrow \infty
}\frac{\left\vert A_{T}\right\vert }{\left\vert \Pi
_{T}\right\vert }, \eea provided the limit exists. We say a random
variable $\xi$ has geometric distribution with parameter $\frac
{1}{2}$ (for brevity, \textbf{geometrically distributed}) if
$\mathbb{P} (\xi = n) = \frac {1}{2^{n}} $ for $n = 1, 2, \dots $.
A consequence of the structure theorem is that
\begin{equation}
\mathbb{P}\left( {x : \gamma _{m}(x) = (k_{1},\dots ,k_{m}) }
\right) \ =\ \frac{1} {2^{k_{1}+\cdots+k_{m}}}\ =\
\prod_{i=1}^{m}\frac{1}{2^{k_{i}}}.
\end{equation} Both the expectation and variance
of a geometrically distributed random variable is $2$. For a seed
$x_{0}$ let $x_{m} = M^{m} (x_{0})$ be the $m$\textsuperscript{th}
iterate. A natural quantity to investigate is
$\frac{x_{m}}{\left(\frac{3}{4}\right)^m x_0}$, where
$\left(\frac{3}{4}\right)^m x_0$ is the expected value of $x_m$.

\begin{theorem}[\cite{KonSi}]\label{thm:xidistrib} The $k$-values are
independent geometrically distributed random variables. Further,
for any $a \in \mathbb{R}$ \be \mathbb{P}\left( \frac{\log_2
\left[ \frac{x_{m}}{\left(\frac{3}{4}\right)^m
x_0}\right]}{\sqrt{2m}} \leq a \right) \ =\ \mathbb{P}\left(
\frac{S_{m}-2m}{\sqrt{2m}} \leq a \right),\ee where $S_{m}$ is the
sum of $m$ geometrically distributed (with parameter $\foh$)
i.i.d.r.v. By the Central Limit Theorem, the right hand side
converges to a Gaussian integral as $m \to \infty$. The paths are
also independent, and so the $\left( 3x+1\right)$-Paths are those
of a geometric Brownian motion with drift $\log \frac{3}{4}$.
\end{theorem}

We remind the reader that a Brownian motion (and hence a geometric
Brownian motion) can be realized as the limit of a random walk;
the same phenomenon occurs here. The drift corresponds to the fact
that the expected value is $\left(\frac{3}{4}\right)^m x_0$,
rather than just $x_0$.

It is worth remarking that a consequence of the drift being
$\log_2 \frac{3}{4}$ (which is negative) is that it is natural to
expect that typical trajectories return to the origin. This
statement extends completely to $(d,g,h)$-Maps discussed in
Appendix B. Theorem \ref{thm:xidistrib} is immediately applicable
to investigations in base two (which is uninteresting as all first
digits are 1). To study the $3x+1$ Problem in base $B$, one simply
multiplies by $\frac{1}{\log_2 B}$, as $\frac{\log_2 x}{\log_2 B}
= \log_2 B$. This replaces $S_m - 2m$ with $\frac{S_m - 2m}{\log_2
B}$ or $(S_m - 2m) \log_B 2$.

\subsection{A Tale of Two Limits}

The $\left( 3x+1\right) $-system, $\overrightarrow{X_{T}}=\left\{
x_{i}\right\} _{0\leq i\leq T}$, is probably not Benford for \textit{any }%
starting seed $x_{0}$ as we expect all of the terms to eventually
be 1. If we stop the sequence after hitting 1 and consider the
proportion of terms having a given leading digit $j$, this is a
rational number, whereas $\log _{10}j$ is not. Of course, this
rational number should be close to $\log_{10}j$, but it is
difficult to quantify this proximity since it is easy to find
arbitrarily large numbers decaying to 1 after even one iteration
of the $\left( 3x+1\right)$-map.

One sense in which Benford behavior can be proved is the same as
the sense in which $\left( 3x+1\right)$-paths are those of a
geometric Brownian motion. We use the structure theorem to prove

\begin{theorem}\label{thm:3x+1isBenford} Let $B$ be any real
number such that $\log_B 2$ is irrational of type $\kappa <
\infty$; for example, one may take any integer $B$ which is not a
perfect power of $2$ (see \eqref{eq:irrtypedef} for a definition
of type $\kappa$ and Theorem \ref{thm:irrtypelog} for a proof of
the irrationality type of such integers). Then for any $[a,b]
\subset [0,1]$, \be
\lim\limits_{m\rightarrow \infty }%
\mathbb{P}\left( \log_B \left[
\frac{x_{m}}{\left(\frac{3}{4}\right)^m x_0}\right] \ \bmod\ 1\in
[a,b] \right) \ = \ b-a.\ee As $\left(\frac{3}{4}\right)^m x_0$ is
the expected value of $x_m$, this implies the distribution of the
ratio of the actual versus predicted value after $m$ iterates
obeys Benford's Law (base $B$). If $B = 2^n$ for some integer $n$,
in the limit $\log_B \left[
\frac{x_{m}}{\left(\frac{3}{4}\right)^m x_0}\right] \ \bmod\ 1$
takes on the values $0, \frac{1}{n}, \frac{2}{n}, \dots,
\frac{n-1}{n}$ with equal probability, leading to a non-Benford
digit bias depending only on $n$.
\end{theorem}

Notice that since probability is defined through density, this is
really two highly non-interchangeable limits: \begin{align} &
\lim_{m\rightarrow \infty }\mathbb{P}\left( \log_B \left[
\frac{x_{m}}{\left(\frac{3}{4}\right)^m x_0}\right]\bmod 1\in
[a,b] \right)\nonumber\\ & \ \  = \ \lim_{m\rightarrow \infty
}\lim\limits_{T\rightarrow \infty }\frac{\#\left\{ x_{0}\in \Pi
_{T}:\log_B \left[ \frac{x_{m}}{\left(\frac{3}{4}\right)^m
x_0}\right]\bmod 1\in [a,b] \right\} }{\#\Pi _{T}}. \end{align}
Though this is completely natural, it is worth remarking for the
sake of precision. Of course, a good starting seed (one with a
long life-span) should give a close approximation of Benford
behavior, just as it will also be a generic Brownian sample path;
this is supported by numerical investigations (see
\S\ref{subsec:numinv3x+1}).

Let $\xi_{1},\xi_{2},\dots $ be independent geometrically
distributed random variables with $\mathbb{P}\left(
\xi_{i}=n\right) =\frac{1}{2^{n}}$, $n=1,2,\dots$, and
$\mathbb{E}\left( \xi_{i}\right) =2$, $\mbox{Var}\left( \xi_{i}\right) =2$. Let $%
S_{m}=\sum_{i=1}^{m}\xi_{i}$. Let $\zeta_{i}=\xi _{i}-2$,
$\overline{S}_{m}=\sum_{i=1}^{m}\zeta _{i}=S_{m}-2m$. We know the
distribution of $\log_B \left[
\frac{x_{m}}{\left(\frac{3}{4}\right)^m x_0}\right]$ is the same
as that of $(S_m - 2m)\log_B 2 = \osm \log_B 2$. The proof is
complicated by the fact that the sum of $m$ geometrically
distributed random variables itself has a binomial distribution,
supported on the integers. This gives a lattice distribution for
which we cannot obtain sufficient bounds on the error, even by
performing an Edgeworth expansion and estimating the rate of
convergence in the Central Limit Theorem. The problem is that the
error in missing a lattice point is of size $\frac1{\sqrt{m}}$,
and we need to sum $\sqrt{m}h(m)$ terms (for some $h(m) \to
\infty$). We are able to surmount this obstacle by an error
analysis of the rate of convergence to equidistribution of $k
\log_B 2 \bmod 1$.

\subsection{Proof of Theorem
\ref{thm:3x+1isBenford}}\label{subsec:proofbenford3x+1}

To prove Theorem \ref{thm:3x+1isBenford} we first collect some
needed results. The proof is similar in spirit to Theorem
\ref{thm:main}, with the needed results playing a similar role as
the three conditions; however, the discreteness of the $3x+1$
problem leads to some interesting technical complications, and it
is easier to give a similar but independent proof than to adjust
notation and show Conditions \ref{cond:1} through \ref{cond:4} are
satisfied.

In the statements below, $[a,b]$ is an arbitrary sub-interval of
$[0,1]$. By the Central Limit Theorem, the distribution of $\osm$
(although it only takes integer values) is approximately a
Gaussian with standard deviation of size $\sqrt{m}$. Let $c \in
\left(0,\frac12\right)$ and set $M = m^c$. Let \be I_\ell\ =\
\{\ell M, \ell M + 1, \dots, (\ell+1)M - 1\}\ee and $C = \log_B 2$
be an irrational number of type $\kappa$ (see
\eqref{eq:irrtypedef}). Soundararajan informed us that one does
not need $\log_B 2$ to be of finite type for our applications. For
integer $B$, if $B^p - 2^q > 0$ then it is at least $1$, and one
obtains $o(M)$ instead of $O(M^\delta)$ in \eqref{eq:EDerror}; the
advantage of using finite type is we obtain sharper estimates on
the rate of convergence, as well as being able to handle
non-integral bases $B$.

Let $\eta(x)$ denote the density of the standard normal: \be
\eta(x)\ =\ \frac1{\sqrt{2\pi}}\ e^{-x^2/2}.\ee We collect some
results needed for the proof of Theorem \ref{thm:3x+1isBenford}:

\begin{itemize}

\item From the Central Limit Theorem (see \cite{Fe}, Chapter XV):
For any $k\in\Z$, \bea\label{eq:ApproxSm} \text{Prob}(C \cdot \osm
= C \cdot k) & \ = \ & \text{Prob}\left(\frac{\osm}{\sqrt{m}} =
\frac{k}{\sqrt{m}}\right) \nonumber\\ & \ = \ &
\frac{1}{\sqrt{m}}\ \eta\left(\frac{k}{\sqrt{m}}\right) +
o\left(\frac{1}{\sqrt{m}}\right). \eea We may write
$o\left(\frac{1}{\sqrt{m}}\right)$ as $O\left(\frac1{\sqrt{m}
g(m)}\right)$ for some monotone increasing $g(m)$ which tends to
infinity. We use this to approximate the probability of $\osm =
k$. For future use, choose any monotone $h(m)$ tending to infinity
such that $h(m)= o\left(g(m)\right)$, $h(m) =
o\left(m^{1/2005}\right)$ and $\frac{h(m)M}{\sqrt{m}} =
o\left(m^{-1/2005}\right)$. As $M = m^c$ with $c < \frac12$, if
$c$ is sufficiently small then such an $h$ exists.

\item Let $k_1, k_2 \in I_\ell$. Then
\begin{align}\label{eq:GaussianChange} & \left| \oosm\
\feta{\frac{k_1}{\sqrt{m}}} - \oosm\ \feta{\frac{k_2}{\sqrt{m}}}
\right|\nonumber\\ & \ \ \ \le \ \oosm\ e^{-\ell^2M^2/2m} \cdot
\left(1 - \exp\left(-\frac{2\ell M^2 + M^2}{2m}\right) \right).
\end{align} In practice this implies that for the $\ell$ we must study,
there is negligible variation in the Gaussian  for $k \in I_\ell$.

\item By Poisson Summation (see page 63 of \cite{Da}),
\be\label{eq:PoissonSummBenford} \frac1{\sigma}
\sum_{n=-\infty}^\infty e^{-n^2\pi/\sigma^2} \ = \
\sum_{n=-\infty}^\infty e^{-n^2 \pi \sigma^2}, \ \ \ \ \sigma > 0.
\ee We often take $\sigma^2 = \frac{2m}{\pi M^2}$, and use this to
calculate the main term (as $\sigma \to \infty$, both sides of
\eqref{eq:PoissonSummBenford} tend to $1$).

\item For any $\epsilon > 0$, letting $\delta = 1 + \epsilon -
\frac1{\kappa} < 1$ we have \be\label{eq:EDerror} \#\{k \in
I_\ell: kC \bmod 1 \in [a,b]\} \ = \ M(b-a) + O(M^{\delta}).\ee
The quantification of the equidistribution of $kC \bmod 1$ is the
key ingredient in proving Benford behavior base $B$ (with $C =
\log_B 2$). The rate of equidistribution, given the finiteness of
the irrationality type of $C$, follows from the Erd\"{o}s-Turan
Theorem. As this is the key argument in our analysis, we provide a
sketch of the proof in Appendix \ref{sec:irrtypelog}; see Theorem
3.3 on page 124 of \cite{KN} for complete details (while the proof
given only applies for $I_0$, a trivial translation yields the
claim for any $I_\ell$).

\end{itemize}

\begin{proof}[Proof of Theorem \ref{thm:3x+1isBenford}]
We must show that as $m\to \infty$, for any $[a,b] \subset [0,1]$,
\be P_m(a,b) \ = \ \text{Prob}(C\osm \bmod 1 \in [a,b]) \ee tends
to $b-a$.  We have \bea\label{eq:Pmabtwopieces} P_{m}(a,b) & \ = \
& \sum_{|\ell| \le \frac{\sqrt{m}h(m)}{M}}
\text{Prob}(\osm = k \in I_\ell: kC \bmod 1 \in [a,b]) \nonumber\\
& & \ + \ \sum_{|\ell| > \frac{\sqrt{m}h(m)}{M}} \text{Prob}(\osm
= k \in I_\ell: kC \bmod 1 \in [a,b]). \eea The second sum in
\eqref{eq:Pmabtwopieces} is bounded by \be\label{PmabsecondsumCLT}
\text{Prob}\left(\osm = k: |k| \ge \frac{\sqrt{m}h(m)}{M}\right).
\ee By the Central Limit Theorem, \eqref{PmabsecondsumCLT} is
$o(1)$. Alternatively, using the techniques below (with $[a,b] =
[0,1]$), one can show $\text{Prob}\left(|\osm| \le
\frac{\sqrt{m}h(m)}{M}\right) = 1 + o(1)$, which implies
\eqref{PmabsecondsumCLT} is $o(1)$. As we are not summing
\eqref{PmabsecondsumCLT}, it is okay to have an error here of size
$\oosm$ (and errors of approximately this size arise if we add or
subtract a lattice point). Therefore \bea P_{m}(a,b) & \ = \ &
\sum_{|\ell| \le \frac{\sqrt{m}h(m)}{M}} \text{Prob}(\osm = k \in
I_\ell: kC \bmod 1 \in [a,b]) + o(1) \nonumber\\ & = &
\sum_{|\ell| \le \frac{\sqrt{m}h(m)}{M}} P_{m,\ell}(a,b) + o(1).
\eea

The proof is completed by showing the above is $b-a + o(1)$.
Consider an interval $I_\ell$. By \eqref{eq:EDerror}, the number
of $k\in I_\ell$ such that $kC \bmod 1 \in [a,b]$ is $(b-a)M +
O(M^\delta)$, $\delta < 1$. By \eqref{eq:ApproxSm}, the
probability of each such $k$ is $\oosm\ \feta{\frac{k}{\sqrt{m}}}
+ O\left(\frac1{\sqrt{m}g(m)}\right)$. We now use
\eqref{eq:GaussianChange} to bound the error from evaluating all
the $\feta{\frac{k}{\sqrt{m}}}$ at $k = \ell M$ and find
\begin{align}\label{eq:fourpieces} P_{m,\ell}(a,b) & =
\frac{(b-a)M}{\sqrt{m}} \left[\feta{\frac{\ell M}{\sqrt{m}}} +
O\left(e^{-\ell^2M^2/2m}\right)\cdot\left(1
- \exp\left(-\frac{2\ell M^2 + M^2}{2m}\right)\right)\right] \nonumber\\
& \ \ \ \ +\ O\left(M \cdot \frac{1}{\sqrt{m}g(m)}\right) +
O\left(M^\delta \cdot \oosm\ \feta{\frac{\ell
M}{\sqrt{m}}}\right); \end{align} summing over all $|\ell| \le
\frac{\sqrt{m}h(m)}{M}$ gives $P_m(a,b)+o(1)$. This gives four
sums, which we must show are $b-a + o(1)$.

The sums over $|\ell| \le \frac{\sqrt{m}h(m)}{M}$ of the first and
fourth pieces of \eqref{eq:fourpieces} are handled by Poisson
Summation. We have for the first piece that
\begin{align}\label{eq:firstpiecetwoterms} & \sum_{|\ell| \le
\frac{\sqrt{m}h(m)}{M}} \frac{(b-a)M}{\sqrt{m}}\ \feta{\frac{\ell
M}{\sqrt{m}}}\nonumber\\  &\ \ \ \ = \ \sum_{\ell=-\infty}^\infty
\frac{(b-a)M}{\sqrt{m}}\ \feta{\frac{\ell M}{\sqrt{m}}} -
\sum_{|\ell| > \frac{\sqrt{m}h(m)}{M}} \frac{(b-a)M}{\sqrt{m}}\
\feta{\frac{\ell M}{\sqrt{m}}}.\end{align} As $h(m) \to \infty$,
the second sum in \eqref{eq:firstpiecetwoterms} is bounded by \be
\int_{|x| \ge \frac{\sqrt{m}h(m)}{M}} \frac{1}{\sqrt{2\pi m/M^2}}\
e^{-x^2/2(m/M^2)}dx \ = \ \frac1{\sqrt{2\pi}} \int_{|u| \ge h(m)}
e^{-u^2/2}du \ = \ o(1). \ee Using \eqref{eq:PoissonSummBenford}
with $\sigma^2 = \frac{2m}{\pi M^2}$ gives \bea \sum_{|\ell| \le
\frac{\sqrt{m}h(m)}{M}} \frac{(b-a)M}{\sqrt{m}} \feta{\frac{\ell
M}{\sqrt{m}}} & \ = \ & (b-a) \sum_{\ell=-\infty}^\infty
\frac{1}{\sqrt{2\pi m/M^2}} \ e^{- \ell^2 /2 (m/M^2)} + o(1)
\nonumber\\ & = & (b-a) \sum_{\ell=-\infty}^\infty e^{-\ell^2
\cdot 2\pi^2 m/M^2} + o(1) \nonumber\\ & = & b-a + O\left(
\frac{e^{-2\pi^2 m/M^2}}{1 - e^{-2\pi^2 m/M^2}}\right) + o(1) \eea
as the final sum over $\ell \neq 0$ is bounded by a geometric
series and $M = m^c$ with $c<\frac12$. Thus the first piece from
\eqref{eq:fourpieces} gives $b-a + o(1)$.

As the Gaussian is a monotone function (for $x \ge 0$ or $x \le
0$), a similar argument shows the sum over $|\ell| \le
\frac{\sqrt{m}h(m)}{M}$ of the fourth piece of
\eqref{eq:fourpieces} contributes $O(M^{\delta - 1}) + o(1)$. It
is here that we use $C \osm$ is a very special equidistributed
sequence modulo $1$, namely it is of the form $k C \bmod 1$. This
allows us to control the discrepancy (how many $k\in I_\ell$ give
$k C \bmod 1 \in [a,b]$).

We must now sum over $|\ell| \le \frac{\sqrt{m}h(m)}{M}$ the
second and third pieces of \eqref{eq:fourpieces}. For the second
piece, we have \be\label{eq:secondpiecesum} \sum_{|\ell| \le
\frac{\sqrt{m}h(m)}{M}} \frac{M}{\sqrt{m}}\ e^{-\ell^2 M^2/2m}
\left[ 1 - \exp\left(-\frac{2\ell M^2 + M^2}{2m}\right)\right].
\ee As $|\ell| \le \frac{\sqrt{m}h(m)}{M}$ and $M = m^c$ with $c <
\foh$, we have \be \frac{2\ell M^2 + M^2}{2m} \ \ll \
\frac{h(m)M}{\sqrt{m}}. \ee Recall we chose $h(m)$ and $c$ such
that $\frac{h(m)M}{\sqrt{m}} = o\left(m^{-1/2005}\right)$.
Therefore \be 1 - \exp\left(-\frac{2\ell M^2 + M^2}{2m}\right) \
\ll \ m^{-1/2005}. \ee As we chose $h(m)$ such that $h(m) =
o\left(m^{1/2005}\right)$, the sum in \eqref{eq:secondpiecesum} is
 \be \ll \ \frac{\sqrt{m}h(m)}{M} \cdot \frac{M}{\sqrt{m}}
\frac1{m^{1/2005}} \ = \ \frac{h(m)}{m^{1/2005}} \ = \ o(1), \ee
proving the second piece in \eqref{eq:fourpieces} is negligible.

We are left with the sum over $|\ell| \le \frac{\sqrt{m}h(m)}{M}$
of the third piece in \eqref{eq:fourpieces}. Its contribution is
\be O\left( \frac{\sqrt{m}h(m)}{M} \cdot \frac{M}{\sqrt{m}g(m)}
\right)\ = \ O\left(\frac{h(m)}{g(m)}\right) \ = \ o(1). \ee

Collecting the evaluations of the sums of the four pieces in
\eqref{eq:fourpieces}, we see that \be P_m(a,b) \ = \ b - a +
o(1), \ee which completes the proof of Theorem
\ref{thm:3x+1isBenford} if $B \neq 2^n$ (and thus proves Benford
behavior base $10$ because, by Theorem \ref{thm:irrtypelog},
$\log_{10}2$ has finite irrationality type).

Consider now the case when $B = 2^n$. As $S_m$ takes on integer
values, the possible values modulo 1 for $(S_m-2m) \log_B 2$ are
$\{0,\frac1{n},\dots,\frac{n-1}{n}\}$. An identical argument shows
each of these values is equally likely; by determining which
intervals $[\log_B d, \log_B(d+1))$ they lie in, one can determine
the (non-Benford) digit bias in this case. See also
\S\ref{subsec:numinv3x+1}. \end{proof}

In Appendix \ref{sec:dghmaps} a generalization of the $3x+1$ map
is discussed; for such systems, one can easily prove the analogue
of Theorem \ref{thm:3x+1isBenford}.

\subsection{Numerical Investigations}\label{subsec:numinv3x+1}

Theorem \ref{thm:3x+1isBenford} implies that the first digit of
$\frac{x_m}{\left(\frac{3}{4}\right)^m x_0}$ will not be Benford
in a base $B = 2^n$. As $S_m$ takes on integer values, $(S_m-2m)
\log_B 2$ is equally likely to be any of
$0,\frac1{n},\dots,\frac{n-1}{n}$. We considered $100,000$ seeds
congruent to 1 modulo 6, starting at $419,753,999,998,525$. We can
rapidly analyze the behavior of such large numbers by representing
each number as an array and then performing the required
operations (multiplication by $3$, addition by $1$ and division by
$2$) digit by digit. Taking $m= 10$, we analyzed the first digits
for $B = 4, 8$ and $16$. We
have (theoretical predictions in parentheses)\\

\begin{center}
\begin{tabular}{|l||r|r|r|r|r|r|r|}
  \hline
  \text{First} &  &  &  &  &  &  &   \\
\text{Digit} & 1 & 2 & 3 & 4 & 5 & 6 & 7  \\
  \hline
\text{Base $4$} & 50.2\% (50.0\%)  & 49.8\% (50.0\%)  & 0\%  & N/A  & N/A  & N/A & N/A   \\
  \hline
  \text{Base $8$} & 33.1\% (33.3\%)  & 33.6\% (33.3\%)  & 0\%  & 33.3\% (33.3\%) & 0\%
   & 0\% & 0\%  \\
  \hline
\end{tabular}
\end{center}

\bigskip

In base 16, we only observe digits 1, 2, 4 and 8; all should occur
$25\%$ of the time; we observe them with frequencies $25.0\%,
25.0\%, 25.3\%$ and $24.8\%$. In base 10, we observe\\

\begin{center}
\begin{tabular}{|l|r|r|r|r|r|r|r|r|r|}
  \hline
  \text{First Digit} & 1 & 2 & 3 & 4 & 5 & 6 & 7 & 8 & 9 \\
  \hline
  \text{Observed} & 29.8\% & 17.9\% & 12.1\% & 10.0\% & 8.5\% & 9.8\% & 2.4\% & 8.7\% & 0.9\%   \\
  \text{Benford} & 30.1\% & 17.6\%  & 12.5\% & 9.7\% & 7.9\% & 6.7\% & 5.8\% & 5.1\% & 4.6\%  \\
  \hline
\end{tabular}
\end{center}

\bigskip

The difficulty in performing these experiments is that our theory
is that of two limits, $T \to \infty$ and then $m\to\infty$. We
want to choose large seeds $x_0$ (at least large enough so that
after $m$ applications of the $3x+1$ map we haven't hit $1$);
however, that requires us to examine (at least on a log scale) a
large number of $x_0$. Taking larger starting values (say of the
order $10^{100}$) makes it impractical to study enough consecutive
seeds. In these cases, to approximate the limit as $T\to\infty$ it
is best to choose $100,000$ seeds from a variety of starting
values and average.

While we cannot yet prove that the iterates of a generic fixed
seed are Benford, we expect this to be so. The table below records
the percent of first digits equal to $j$ base $10$ for a 100,000
random digit number under the $3x+1$ map (as the $3x+1$ map
involves simple digit operations, we may represent numbers as
arrays, and the computations are quite fast). We performed two
experiments: in the first we removed the highest power of $2$ in
each iteration ($799,992$ iterates), while in the second we had
$M(x) = 3x+1$ for $x$ odd and $\frac{x}{2}$ for $x$ even
($2,402,282$ iterates). In both, the observed probabilities are
extremely close to the Benford predictions (for each digit, the
corresponding $z$-statistics range from about $-2$ to $2$).


\begin{center}
\begin{tabular}{|l|r||r|r||r|r|}
  \hline
First   &   Benford &       &  &  &
  \\
Digit   &   Probability &   Removing 2    & $z$-statistic & Not
Removing $2$  &
  $z$-statistic  \\

\hline
1   &   0.3010    &   0.3021     &   2.00    &   0.3012    &   0.63    \\
2   &   0.1761    &   0.1752    &   -2.10   &   0.1763    &   0.98    \\
3   &   0.1249    &   0.1242    &   -1.97   &   0.1248    &   -0.69   \\
4   &   0.0969    &   0.0967   &   -0.50   &   0.0967   &   -1.14   \\
5   &   0.0792    &   0.0792   &   0.03    &   0.0792   &   -0.06   \\
6   &   0.0670    &   0.0671   &   0.56    &   0.0667   &   -1.32   \\
7   &   0.0580    &   0.0582   &   0.68    &   0.0581    &   0.89    \\
8   &   0.0512    &   0.0513   &   0.79    &   0.0510   &   -0.77   \\
9   &   0.0458    &   0.0460   &   0.99    &   0.0459   &   1.02    \\
  \hline
\end{tabular}
\end{center}

\bigskip

We calculated the $\chi^2$ values for both experiments: it is
$12.38$ in the first ($M(x) = \frac{3x+1}{2^k}$) and $6.60$ in the
second ($M(x) = 3x+1$ for $x$ odd and $\frac{x}2$ otherwise). As
for $8$ degrees of freedom, $\ga = .05$ corresponds to a $\chi^2$
value of $15.51$, and $\ga = .01$ corresponds to $20.09$, we do
not reject the null hypothesis and our experiments support the
claim that the iterates of both maps obey Benford's law.

\section{Conclusion and Future Work}

The idea of using Poisson Summation to show certain systems are
Benford is not new (see for example \cite{Pin} or page 63 of
\cite{Fe}); the difficulty is in bounding the error terms. Our
purpose here is to codify a certain natural set of conditions
where the Poisson Summation can be executed, and show that
interesting systems do satisfy these conditions; a natural future
project is to determine additional systems that can be so
analyzed. One of the original goals of the project was to prove
that the first digits of the terms $x_m$ in the $3x+1$ Problem are
Benford. While the techniques of this paper are close to handling
this, the structure theorem at our disposal makes
$\frac{x_m}{\left(\frac34\right)^mx_0}$ the natural quantity to
investigate (although numerical investigations strongly support
the claim that for any generic seed, the iterates of the $3x+1$
map are Benford); however, we have not fully exploited the
structure theorem and the geometric Brownian motion, and hope to
return to analyzing the first digit of $x_m$ at a later time.
Similarly, additional analysis of the error terms in the
expansions and integrations of $L$-functions may lead to proving
Benford behavior on the critical line, and not just near it,
although our results on values of $L$-functions near the critical
line as well as the digits of values of characteristic polynomials
of random matrix ensembles support the conjectured Benford
behavior.

\section*{Acknowledgements}

We thank Arno Berger, Ted Hill, Ioannis Karatzas, Jeff Lagarias,
James Mailhot, Jeff Miller, Michael Rosen, Yakov Sinai and Kannan
Soundararajan for many enlightening conversations, Dean Eiger and
Stewart Minteer for running some of the $3x+1$ calculations, Klaus
Schuerger for pointing out some typos in an earlier draft, and the
referee for many valuable comments (especially suggesting we study
characteristic polynomials of unitary matrices). Both authors
would also like to thank the 2003 Hawaii International Conference
on Statistics and The Ohio State University for their hospitality,
where much of this work was written up.


\appendix

\section{Values of Characteristic
Polynomials}\label{sec:appvalcharpoly}

Consider the random matrix ensemble of $N\times N$ unitary
matrices $U$ (with eigenvalues $e^{\i\theta_n}$) with respect to
Haar measure; the probability density of $U$ is \be p_N(U) \ = \
\frac1{(2\pi)^NN!}\ \prod_{1\le j<m \le N}\
\left|e^{\i\theta_j}-e^{\i\theta_m}\right|. \ee  Let \be
Z(U,\theta) \ = \ \det(I - Ue^{-\i\theta}) \ = \ \prod_{n=1}^N
\left(1 - e^{\i(\theta_n-\theta)}\right) \ee be the characteristic
polynomial of $U$. The values of characteristic polynomials have
been shown to be a good model for the values of $L$-functions. Of
interest to us are the results in \cite{KeSn}, where an analogue
of the log-normal law of $L$-functions (Theorem
\ref{thm:HejLogNorm}) is shown for random matrix ensembles: as
$N\to\infty$ the average of the absolute value of the
characteristic polynomials of unitary matrices is Gaussian.
Specifically, let $\rho_N(x)$ be the probability density for $\log
|Z(U,\theta)|$ averaged with respect to Haar measure (Equation
(36) of \cite{KeSn}), and set \be\label{eq:rhoNwtrhoN}
\widetilde{\rho}_N(x) \ = \ \sqrt{Q_2(N)}\ \rho_N(\sqrt{Q_2(N)}\
x). \ee Here $Q_2(N)$ is the variance, and by Equation (11) of
\cite{KeSn} satisfies \be Q_2(N) \ = \ \frac{\log N}2 +
\frac{\gamma+1}2 + \frac1{24N^2} + O(N^{-4}). \ee Equation (53) of
\cite{KeSn} (and the comment immediately after it) yield

\begin{theorem}[Keating-Snaith] With $\widetilde{\rho}_{N}$ as above, \be
\widetilde{\rho}_N(x)dx \ = \ \frac1{\sqrt{2\pi}}\ e^{-x^2/2}dx +
O\left( (\log N)^{-3/2}dx\right). \ee
\end{theorem} In terms of $\rho_N$, from \eqref{eq:rhoNwtrhoN}
we immediately deduce that \be \rho_N(x)dx \ = \ \frac1{\sqrt{2\pi
Q_2(N)}}\ e^{-x^2 / 2 Q_2(N)}dx + O\left( Q_2(N)^{-2}dx\right);
\ee note the pointwise errors are of size one over the square of
the variance. It is easy to show the conditions of Theorem
\ref{thm:poissum} are satisfied. These errors are significantly
smaller than the number theory analogues, in part due to the
additional averaging (the formulas here are for averages with
respect to Haar measure, whereas in number theory we studied one
specific $L$-function). We thus have

\begin{theorem}\label{thm:valuesCharPoly}
As $N\to\infty$, the distribution of digits of the absolute values
of the characteristic polynomials of $N\times N$ unitary matrices
(with respect to Haar measure) converges to the Benford
probabilities.
\end{theorem}

\begin{proof} As the main term is given by a Gaussian, the only difficulty
is in verifying Conditions \ref{cond:1} and \ref{cond:4}. In our
current setting, $\sqrt{Q_2(N)}$ is playing the role of $T$. Let
$h(N) = \log Q_2(N)$. As \be
\int_{-\sqrt{Q_2(N)}h(N)}^{\sqrt{Q_2(N)}h(N)} \frac1{\sqrt{2\pi
Q_2(N)}}\ e^{-x^2 / 2 Q_2(N)}dx  \ = \ 1 + o(1), \ee Condition
\ref{cond:1} is satisfied. For Condition \ref{cond:4}, note
$E_T(b+k)-E_t(a+k)$ becomes $O\left(Q_2(N)^{-2}\right)$, and thus
\be \sum_{|k| \le \sqrt{Q_2(N)} h(N)} \left[E_N(b+k) -
E_N(a+k)\right] \ \ll \ \sqrt{Q_2(N)} h(N) Q_2(N)^{-2} \ \ll \
\frac{\log Q_2(N)}{Q_2(N)^{3/2}}. \ee
\end{proof}

\begin{remark} While we believe the distribution of digits of
$L$-functions on the critical line is Benford, our results
(Theorem \ref{thm:benfordLfns} and Corollary
\ref{cor:benfordLfns}) apply to values just off the critical line.
Theorem \ref{thm:valuesCharPoly} may thus be interpreted as
providing additional support to the conjectured Benford behavior
of $L$-functions on the critical line.
\end{remark}

\begin{remark} In our earlier investigations of Benford behavior,
we used either the counting measure (first $N$ terms of a
sequence) or Lebesgue measure (values of the function at arguments
$t\in [0,N]$), with $N\to\infty$. We have an extra averaging here.
We are not looking at the characteristic polynomials of a sequence
of unitary matrices $U_N$ (where $U_N$ is $N\times N$). Instead
for each $N$ we use Haar measure on $N\times N$ unitary matrices
to average the values of the characteristic polynomials, and then
send $N\to\infty$. The averaged characteristic polynomials play an
analogous role to our $L$-functions from before.
\end{remark}


\section{Irrationality type of $\log_B 2$ and Equidistribution}\label{sec:irrtypelog}

\begin{theorem}\label{thm:irrtypelog} Let $B$ be a positive integer
not of the form $2^n$ for an integer $n$. Then $\log_B 2$ is of
finite type. \end{theorem}

\begin{proof} By \eqref{eq:irrtypedef},
we must show for some finite $\kappa > 0$ that \be \left| \log_B 2
- \frac{p}{q}\right| \ \gg \ \frac1{q^\kappa}. \ee As \be \left|
\frac{\log 2}{\log B} - \frac{p}{q}\right|\ =\ \frac{|q \log 2 - p
\log B|}{|q| \log B}, \ee it suffices to show $|q\log 2 - p \log
B| \gg q^{-\kappa'}$. This follows immediately from Theorem 2 of
\cite{Ba}, which implies that if $\ga_j$ and $\beta_j$ are
algebraic integers of heights at most $A_j (\ge 4)$ and $B (\ge
4)$, then if $\Lambda = \beta_1 \log \ga_1 + \cdots + \beta_n \log
\ga_n \neq 0$, $|\Lambda| > B^{-C\Omega \log \Omega'}$, where $d$
is the degree of the extension of $\Q$ generated by the $\ga_j$
and $\beta_j$, $C = (16nd)^{200n}$, $\Omega = \log \ga_1 \cdots
\log \ga_n$ and $\Omega' = \Omega / \log \ga_n$. We take $B$ to be
maximum of $\beta_1 = q$ and $\beta_2 = -p$. (As stated we need
$\ga_1, \ga_2 \ge 4$; we replace $q \log 2 - p \log B$ with $\foh(
q \log 4 - p \log B^2)$). In our case $d = 1$, $n=2$, $\ga_1 = 4,
\ga_2 = B^2$.  As $B$ is not a power of $2$, $q\log 4 - p \log B^2
\neq 0$ unless $p, q = 0$. In particular, \be\label{eq:appdiscDN}
\left| \log_{B} 2 - \frac{p}{q}\right| \ \gg \
\frac{1}{q^{1+C\Omega\log \Omega'}}. \ee For $B = 10$ we may take
$\kappa = 2.3942 \times 10^{602}$ (though almost surely a lower
number would suffice).
\end{proof}

We show the connection between the irrationality type of $\ga$ and
equidistribution of $n\alpha \bmod 1$; see Theorem 3.3 on page 124
of \cite{KN} for complete details. Define the discrepancy of a
sequence $x_n$ ($n \le N$) by \be D_N \ = \ \frac1{N}\sup_{[a,b]
\subset [0,1]} \left| N(b-a) - \#\{n\le N: x_n \bmod 1 \in
[a,b]\}\right|. \ee The Erd\"{o}s-Turan Theorem (see \cite{KN},
page 112) states that there exists a $C$ such that for all $m$,
\be D_N \ \le \ C\left( \frac1{m} + \sum_{h=1}^m \frac1{h} \left|
\frac1{N} \sum_{n=1}^N e^{2\pi \i h x_n}\right| \right). \ee If
$x_n = n\alpha$, then the sum on $n$ above is bounded by
$\min\left(N,\frac1{|\sin \pi h \alpha|}\right)$ $\le$
$\min\left(N,\frac1{2||h\alpha||}\right)$, where $||x||$ is the
distance from $x$ to the nearest integer. If $\ga$ is of finite
type, this leads to $\sum_{h=1}^m \frac1{h||h\alpha||}$. For $\ga$
of type $\kappa$, this sum is of size $m^{\kappa - 1 +\gep}$, and
the claimed equidistribution rate follows from taking $m = \lfloor
N^{1/\kappa}\rfloor$.

\section{$(d,g,h)$-Maps}\label{sec:dghmaps}
The Benford behavior of $3x+1$ also occurs in $(d,g,h)$-Maps,
defined as follows. Consider positive coprime integers $d$ and
$g$, with $g>d\geq 2$, and a periodic function $h\left( x\right)$
satisfying:

\begin{enumerate}
\item  $h\left( x+d\right) =h\left( x\right) $,

\item  $x+h\left( x\right) \equiv 0  \bmod d$,

\item  $0<\left| h\left( x\right) \right| <g$.
\end{enumerate}

The map $M$ is defined by the formula \be M\left( x\right) \ = \
\frac{gx+h\left( gx\right) }{d^{k}}, \ee where $k$ is uniquely
chosen so that the result is not divisible by $d$. Property (2) of
$h$ guarantees $k\geq 1$. The natural domain of this map is the
set $\Pi $ of positive integers not divisible by $d$ and $g$. Let
$E$ be the set of integers between $1$ and $dg$ that divide
neither $d$ nor $g$, so we can write $\Pi =dg\Bbb{Z}^{+}+E$. The
size of $E$ can easily be calculated: $\left| E\right| =\left(
d-1\right) \left( g-1\right) $. In the same way as before, we have
$m$-paths, which are the values of $k$ that appear in iterations
of $M$, and we again denote them by $\gamma _{m}\left( x\right) $.

The $3x+1$ Problem corresponds to $g=3$, $d=2$, and $h\left(
1\right) =1$, the $3x-1$ Problem corresponds to $g=3$, $d=2$, and
$h\left(
1\right) =-1$, the $5x+1$ Problem corresponds to $g=5$, $d=2$%
, and $h\left( 1\right) =1$, and so on. Similar to Theorem
\ref{thm:xidistrib}, one can show

\begin{theorem}[\cite{KonSi}]
The $(d,g,h)$-Paths are those of a geometric Brownian motion with
drift $\log g - \frac{d}{d-1} \log d$.
\end{theorem}

We expect paths to decay for negative drift and escape to infinity
for positive drift. All results on Benford's Law for the
$(3x+1)$-Problem, in particular Theorem \ref{thm:3x+1isBenford},
generalize trivially to all $(d,g,h)$-Maps, with the
(irrationality) type of $\log_B d$ the generalization of the
(irrationality) type of $\log_B 2$; note Theorem
\ref{thm:irrtypelog} is easily modified to analyze $\log_B d$.

\bigskip

\end{document}